\documentclass[12pt,a4paper]{article}

\usepackage{latexsym}
\usepackage{amsfonts}
\usepackage{amssymb}

\newcommand{\be}{\begin{equation}}
\newcommand{\ee}{\end{equation}}
\newcommand{\bea}{\begin{eqnarray}}
\newcommand{\eea}{\end{eqnarray}}

\renewcommand{\theequation}{\arabic{section}.\arabic{equation}}

\def\G{{\cal G}}                                 %
\def\F{{\cal F}}                                 %
\def\D{{\cal D}}                                 %
\def\I{{\cal I}}                                 %
\def\A{{\cal A}}                                 %
\def\B{{\cal B}}                                 %
\def\M{{\cal M}}                                 %
\def\K{{\cal K}}                                 %
\def\mi{{\mathrm{i}}}                            %
\def\R{{\mathbf R}}                              %
\def\bbR{{\mathbb R}}                            %
\def\cR{{\cal R}}                                %
\def\cL{{\cal L}}                                %
\def\ad{{\mathrm{ad}}}                           %
\def\Ri{{R^{\mathrm{i}}}}                        %
\def\Rit{{R^{\mathrm{i}\theta}}}                 %
\def\CG{{\check G}}                              %
\def\half {\textstyle{\frac{1}{2}}}              %
\def\quarter{\textstyle{\frac{1}{4}}}            %
\def\dt{\left.\frac{d}{dt}\right\vert_{t=0}}     %
\begin{document}

\vspace*{0.5cm}
\begin{center}
{\Large \bf The non-Abelian momentum map for Poisson-Lie symmetries
        on the chiral WZNW phase space}
\end{center}

\vspace{1.0cm}

\begin{center}
L. Feh\'er${}^{a,}$\footnote{Postal address: MTA KFKI RMKI,
H-1525 Budapest 114, P.O.B. 49,  Hungary}
 and I. Marshall${}^b$ \\

\bigskip

{\em
${}^a$Department of Theoretical Physics\\
 MTA  KFKI RMKI and University of Szeged \\
E-mail: lfeher@rmki.kfki.hu\\

\bigskip

${}^b$
Department of Mathematics, EPFL \\
1015 Lausanne, Switzerland\\
E-mail: ian.marshall@epfl.ch
}
\end{center}

\vspace{1.0cm}

\begin{abstract}
The  gauge action of the Lie group $G$ on the chiral WZNW phase space $\M_\CG$
of quasiperiodic fields with $\CG$-valued monodromy, where $\CG\subset G$ is an open 
submanifold,
is known  to be a Poisson-Lie (PL) action with respect to any coboundary PL structure
 on $G$, if
the Poisson bracket on $\M_\CG$ is defined by a suitable monodromy dependent 
exchange $r$-matrix.
We describe the momentum map for these symmetries
when $G$ is either a factorisable PL group or
a compact simple Lie group with its standard  PL structure.
The main result is an explicit
one-to-one correspondence between
the monodromy variable $M \in \CG$ and a conventional 
variable $\Omega \in G^*$.
This permits us to convert the PL groupoid associated
with a WZNW exchange $r$-matrix into a `canonical' PL groupoid
constructed from the Heisenberg double of
$G$, and consequently to obtain a natural PL generalization of the 
classical dynamical Yang-Baxter equation.

\end{abstract}

\bigskip



\newpage

\section{Introduction with a brief review of chiral WZNW}

The Wess-Zumino-Novikov-Witten (WZNW) model of conformal field theory \cite{Wi,CFT} 
not only has many important physical applications, but it is also a  rich source of 
intriguing mathematical structures. 
In particular, the investigations \cite{Fad,BDF,AF,FG,Feld}
of the quantum group and Poisson-Lie (PL) symmetries of the chiral 
sectors of the model 
proved to be instructive for developing the mathematical theory of dynamical 
Yang-Baxter equations \cite{EV,ES,E}.

In the work presented here we concentrate our attention on a result of 
\cite{BFP} in which, given an arbitrary coboundary PL structure on the 
underlying finite dimensional Lie group $G$, there were found PL 
$G$-symmetries on the chiral WZNW phase space $\M_\CG$.
If $G$ is equipped with a Poisson bracket (PB) defined in the
usual way by means of a constant $r$-matrix (called $R^\nu$ below), 
then the PB on $\M_\CG$ can be adjusted in such a way that the 
standard gauge action of $G$ on $\M_\CG$ be a PL action.
The analysis of this phenomenon led to the differential equation 
(\ref{1.12}) below, for which a family of solutions was found.

The main purpose of the present paper is to give an explicit
description of the momentum map 
for the above-mentioned action of $G$ on $\M_\CG$.
This allows us to understand the proper geometric meaning of 
(\ref{1.12}), and thus to obtain a natural PL generalization 
(equation (\ref{PL-CDYBE}))
of the classical dynamical Yang-Baxter equation 
(CDYBE) studied in \cite{EV}.
Interestingly, the PL-CDYBE
(\ref{PL-CDYBE}) also appears in a different context in \cite{DM}.

Our second purpose is to find canonical models for certain finite dimensional
PL groupoids associated with the chiral WZNW phase space in \cite{BFP}.
The study of the latter finite dimensional problem is actually equivalent 
to the former infinite dimensional one.  
The result will be derived by translating the momentum map on the infinite 
dimensional phase space $\M_\CG$   
into a corresponding momentum map on the
associated finite dimensional PL groupoid, which has the same PL symmetry.
One could proceed the other way round, that is to first derive 
the momentum map for the finite dimensional groupoid and then obtain the
corresponding result for $\M_\CG$ as a corollary.  
We have chosen to start from the infinite dimensional case, since we wish    
to directly connect with the analysis in \cite{BFP} that motivated this work.

To explain the content of this paper in more detail,
and to fix notations for later reference,
we next present a brief review of the PL symmetries found in \cite{BFP}. 

The WZNW model as a classical field theory on the cylinder \cite{Wi}
can be defined for any (real or complex) Lie group $G$
whose Lie algebra $\G$ is self-dual in the sense that it is  equipped
with an invariant, symmetric, non-degenerate
bilinear form $\langle\ ,\ \rangle$.
The solution of the field equation for the $G$-valued
WZNW field, which is $2\pi$-periodic in the space variable, is
given by the product of left- and right-moving quasi-periodic chiral WZNW fields.
One obtains the chiral WZNW phase space $\M_\CG$, for which the `monodromy matrix' $M$
is restricted to some open submanifold $\CG\subseteq G$,
\begin{equation}
\M_\CG:= \{  \eta \in C^\infty({\bbR}, G)  \,\vert\,
\eta(x + 2\pi) = \eta(x) M, \quad  M\in \CG\}.
\label{1.1}\end{equation}
The space  $\M_\CG$ is equipped with a weakly non-degenerate symplectic form \cite{FG},  
which allows one to associate Hamiltonian vector fields to certain Hamiltonians. 
It was found in \cite{BFP} (see also \cite{Montreal}) that 
the fundamental `admissible Hamiltonians' are 
`smeared out matrix elements' of the chiral field $\eta$, whose Hamiltonian vector fields   
are encoded in the distribution sense by the following  `Poisson bracket 
relations'\footnote{
Clearly, (\ref{1.2}) should generate a Poisson algebra for certain functions of $\eta$ 
to make ($\M_G, \{\ ,\ \}^r_{WZ}$) into a Poisson manifold in the strict sense.
We do not pursue this issue here, since the specification 
of a Poisson algebra is not needed for the purposes of this paper.
For further explanation, see Remark 1.
Note also that, as described in \cite{BFPinPLA}, by setting $r(M)=0$ 
in (\ref{1.2}) one obtains a
quasi-Poisson structure \cite{AKM} on $\M_G$.}
\begin{equation}
\left\{\eta_1(x) , \eta_2(y)\right\}_{WZ}^r
=\eta_1(x)\eta_2(y)\left(
\half \hat I \,{\mathrm{sign}}\,(y-x) + r(M)  \right), \quad 0< x,y<2\pi.
\label{1.2}\end{equation}
We use the standard St Petersburg notation,
$\eta_1(x)=\eta(x)\otimes 1$, $\eta_2(y)=1\otimes \eta(y)$ and so on, 
together with automatic summation over repeated indices.
The interesting object in (\ref{1.2}) is the `exchange $r$-matrix'
$r(M)= r^{ab}(M) T_a \otimes T_b \in \G\wedge \G$;
$\hat I=T_a \otimes T^a$ with $\{ T_a\}$ and $\{T^a\}$
denoting dual bases of $\G$, $\langle T_a, T^b\rangle =\delta_a^b$.
The monodromy matrix $M$ is a $G$-valued function on $\M_\CG$,
since $M= \eta^{-1}(x) \eta(x+2\pi)$ $\forall x\in \bbR$.
Any smooth function $\psi$ on $G$ gives rise  to an admissible 
Hamiltonian $\hat\psi$ on $\M_\CG$ by $\hat \psi(\eta):= \psi(M)$,
and the corresponding Hamiltonian vector fields are encoded all together by the relation 
\be
\{ \eta_1(x) , M_2\}_{WZ}^r =
\eta_1(x)\left( M_2 r^+(M) - r^-(M) M_2\right).
\label{1.3}\ee
As a consequence of (\ref{1.3}), one also has the Poisson brackets 
\be
\{ M_1 {,} M_2\}_{WZ}^r =
 M_1M_2 r(M) + r(M)M_1 M_2 - M_1 r^-(M) M_2 - M_2 r^+(M) M_1,
\label{1.4}\ee
where  $r^\pm := r \pm \frac{1}{2} \hat I$.
The Jacobi identity condition for
the bracket $\{\ ,\ \}_{WZ}^r$, evaluated for 3 smeared out matrix
elements of $\eta$, is equivalent to the `$G$-CDYBE' for $r$:
\begin{equation}
\left[ r_{12}(M), r_{23}(M)\right]
+T_1^a \left( \half \D_a^+ +
r_a^{{\phantom{a}}b}(M) \D_b^-\right)  r_{23}(M)
+ \hbox{cycl. perm.}=
-\quarter \hat f.
\label{1.5}\end{equation}
Here $\hat f:= f_{a b}^{\phantom{ab}c} T^a \otimes T^b \otimes T_c$ with
$[T_a, T_b]=f_{ab}^{\phantom{ab}c}T_c$,
$r_{23} = r^{ab}  (1\otimes T_a\otimes T_b)$ and
$T^a_1= T^a \otimes 1 \otimes 1$ as usual;
for any  function $\psi$ on $G$ we use
\begin{equation}
\D_a^{\pm} = \cR_a \pm {\cal L}_a,
\quad
({\cal R}_a \psi)(M):= \frac{d}{d t}
 \psi(Me^{t T_a} )\Big\vert_{t=0},
\quad
({\cal L}_a \psi)(M):= \frac{d}{d t}
 \psi(e^{t T_a} M)\Big\vert_{t=0}.
\label{1.6}\end{equation}
Equation (\ref{1.5}) admits a distinguished family of
solutions associated with PL symmetries acting
on the chiral WZNW phase space.

Suppose that $R^\nu\in \G\wedge \G$ is a constant $r$-matrix satisfying
\be
[R^\nu_{12}, R^\nu_{23}] + \hbox{cycl. perm.} = - \nu^2 \hat f,
\label{1.7}\ee
where $\nu$ is a numerical parameter.
Then $G$ is a PL group when equipped with the Sklyanin PB
\be
\left\{q_1, q_2\right\}_G^{R^\nu} = [ q\otimes q, R^\nu].
\label{1.8}\ee
The natural right-action of $G$ on $\M_\CG$  is defined by
\be
\M_{\CG}\times G \ni (\eta , q)\mapsto \eta q \in\M_{\CG}.
\label{1.9}\ee
This yields a PL action, for any Poisson algebra generated by (\ref{1.2}), if the exchange
$r$-matrix has the form
\be
r(M)= R^\nu +K^\nu(M),
\label{1.10}\ee
where $K^\nu$ is subject to the equivariance condition
\be
K^\nu(qMq^{-1}) = (q\otimes q)
K^\nu(M) (q^{-1}\otimes q^{-1}).
\label{1.11}\ee
In this case (\ref{1.5}) gives the `PL-CDYBE' for $K^\nu$, 
\begin{equation}
\left[ K_{12}^\nu(M), K_{23}^\nu(M)\right]
-\half T_1^a  \D_a^+
K_{23}^\nu(M)
+ \hbox{cycl. perm.}=
(\quarter -\nu^2) \hat f.
\label{1.12}\end{equation}
In a neighbourhood of $e\in G$, (\ref{1.11}) can be ensured 
by the ansatz 
\be
K^\nu(M) = \langle T_a, f_\nu(\ad_m) T_b \rangle T^a \otimes T^b,
\qquad
m=\log M,
\label{1.13}\ee
where $f_\nu(z)$ is assumed to be an odd analytic function
in a neighbourhood of zero.
A solution \cite{BFP} to (\ref{1.12}) (subsequently shown  
in \cite{FM} to be unique within the
ansatz
(\ref{1.13})) is provided by the function
\be
f_\nu(z) = z^{-1} \left[ \chi(\half z) -\chi(\nu z) \right]
\quad\hbox{with}\quad
\chi(z)=z \coth z.
\label{1.14}\ee

In the present paper we wish to describe explicitly the
non-Abelian momentum map \cite{Lu} 
(for reviews, see also \cite{FG,BB,SoibAMS})   
 that generates the
action of the PL group $(G, \{\ ,\ \}_G^{R^\nu})$ on $(\M_{\CG}, \{\ ,\ \}_{WZ}^r)$
with $r(M)$ given by (\ref{1.10}), (\ref{1.13}), (\ref{1.14}).
Let $(G^*,\{\ ,\ \}_{G^*}^{R^\nu})$ denote the PL dual of $(G, \{\ ,\ \}_G^{R^\nu})$.
The momentum map we are looking for is required to be a
Poisson map $\Omega:   \M_{\CG} \rightarrow G^*$.
It must generate the infinitesimal version of the action in (\ref{1.9}),
which means that it must satisfy the condition
\be
\left( \{ \eta_{ij}(x), \Omega\}^r_{WZ}
\Omega^{-1}, T \right) = (\eta(x) T)_{ij},
\qquad \forall T\in \G.
\label{1.15}\ee
In this formula $\eta_{ij}(x)$ is an arbitrary matrix element of
$\eta(x)\in G$, $\{ \eta_{ij}(x), \Omega\}^r_{WZ} \Omega^{-1}\in \G^*$
is evaluated on $T\in \G$ in the natural manner;
and for simplicity of writing we pretend that we are dealing with matrix Lie groups.
The left hand side of (\ref{1.15}) encodes the Hamiltonian vector fields generated 
on $\M_\CG$ by the matrix elements of $\Omega$, which must be admissible Hamiltonians. 

The momentum map is already known explicitly in two special cases.
If $\nu =0$ with $R^\nu=0$, then the PL symmetry reduces
to classical $G$-symmetry and $f_\nu$ becomes
$f_0(z)= \frac{1}{2}\coth \frac{z}{2} - \frac{1}{z}$
defining the `canonical' (Alekseev-Meinrenken) 
$r$-matrix by (\ref{1.13}) \cite{EV,AM,BFP,PF}.
Now $G^*$ is the Abelian Lie group $\G^*$, and we can
identify it with $\G$ by means of the scalar product.
In terms of the function $m=\log M$ on $\M_{\CG}$,  the relations (\ref{1.3}),
(\ref{1.4}) can be rewritten  \cite{BFP} as
\begin{equation}
\{ \eta(x), m_a \}_{WZ}^r = \eta(x) T_a,
\qquad
\{ m_a, m_b\}_{WZ}^r =-
f_{ab}^{\phantom{ab}c} m_c
\quad\hbox{with}\quad
m_a=\langle T_a, m\rangle.
\label{1.16}\end{equation}
Thus in this case the `Abelian' momentum map is
given by $m=\log M$.
The second well understood case is that of $\nu=\frac{1}{2}$,
when $r(M)$ equals the constant $r$-matrix  $R^{\frac{1}{2}}$.
In this case (\ref{1.4}) is 
recognized to be the Semenov-Tian-Shansky PB \cite{STS},
which is the PB on $G^*$ by using the standard
identification of $G^*$  with an open submanifold $G$.
Correspondingly, for $\nu=\frac{1}{2}$ the momentum map is directly
furnished by the monodromy matrix \cite{FG}.

The above two examples, together with the fact that the
momentum map is the monodromy matrix for many
integrable systems, lead us to expect that for all cases 
the momentum map for the PL symmetries on
the chiral WZNW phase space be provided by the monodromy matrix\footnote{
Then 
(\ref{1.15}) translates into a finite
dimensional problem for the function $\Omega(M)$ thanks to (\ref{1.3}),
and the Poisson brackets between the matrix elements of $\Omega$ 
are well-defined thanks to (\ref{1.4}).}.
In this paper we confirm this expectation for any $\nu\neq 0$
by explicitly exhibiting the required Poisson map $M\mapsto \Omega(M)$,
where $\Omega$ is a convenient coordinate
on the Poisson space $(G^*,\{\ ,\ \}_{G^*}^{R^\nu})$.
Our main result is the strikingly simple formula
\be
M\mapsto \Omega(M)=M^{2\nu}= \exp(2 \nu \log M),
\label{1.17}\ee
which clearly generalizes the $\nu=\frac{1}{2}$ case.
The precise statement is formulated in Theorem 1 for
the factorisable PL groups
and in Theorem 2 for the 
standard compact PL groups.

In \cite{BFP} a finite dimensional PL groupoid had been associated with every WZNW exchange $r$-matrix.
On the other hand \cite{FM}, a PL groupoid can be constructed from the 
Heisenberg double \cite{STS} of the PL group 
$(G, \{\ ,\ \}_G^{R^\nu})$ by forgetting the relationship between
the left- and right-momenta and shifting the PB of the $G$-valued variable  
by a dynamical $r$-matrix term.
Our next result establishes an isomorphism between these two PL groupoids 
by using (\ref{1.17}) to identify
$M\in \CG$ and $\Omega\in G^*$.
See Proposition 1 for the factorisable case and the discussion around (\ref{3.23})
for the compact case.
By the same change of variables, the PL-CDYBE (\ref{1.12}) can be seen to be  a natural
generalization of the CDYBE on $\G^*$.
This formulation of the PL-CDYBE appears in eqs.~(\ref{PL-CDYBE}) and (\ref{3.27}).

The rest of the paper is organized as follows.
In the next section $G$ is assumed to be a  factorisable PL group,
which means that the constant $r$-matrix in (\ref{1.7})
has the form $R^\nu=2\nu R$ ($\nu\neq 0$) with both $R$ and $R^\nu$ belonging to $\G\wedge\G$.
Such $r$-matrices exist for the complex simple Lie groups and their split real forms
as well as for some other special real forms, but  not for 
the compact one \cite{BD,Rawns}.
In the physically most important case of compact
simple Lie groups, all solutions of
(\ref{1.7}) belong to a purely imaginary parameter $\nu$.
The compact case is studied in Section 3 by taking
$R^\nu = \theta \Ri$, where $\theta$ is real and $\Ri\in \G\wedge \G$
is $\mathrm{i}$-times the Drinfeld-Jimbo $r$-matrix
(normalized so that $\nu = \mathrm{i} \theta$ in (\ref{1.7}), see (\ref{3.6})).
This is not a serious restriction of generality since
in the compact case the most general constant $r$-matrix
with non-zero $\nu$ is obtained \cite{Soib} by adding a purely Cartan piece to
a multiple of $\Ri$.
Both in Sections 2 and 3 we shall proceed by comparing the PBs in
(\ref{1.3}), (\ref{1.4}) to corresponding PBs on the Heisenberg double of the
PL group $G$. This will not only simplify the analysis, but it will also
allow us to directly relate the PL groupoids of \cite{BFP} to the Heisenberg double.
An appendix has been added collecting together all relevant technical
details concerning the Heisenberg double of the standard compact PL
groups.
The material in the appendix is standard and fairly well-known, 
but we thought it useful to include it as it renders our 
description of the compact case essentially self-contained.
Finally, a summary of the results is contained in Section 4.

\noindent
{\em  Remark 1.}
Consider a finite dimensional representation $\Lambda: G\rightarrow GL(V)$ of $G$ and 
a smooth, $2\pi$-periodic  function $\phi: \bbR \rightarrow \mathrm{End}(V)$ 
subject to $\phi^{(k)}(0)=\phi^{(k)}(2\pi)=0$ for every integer $k\geq 0$. 
The `smeared out matrix element' alluded to around (\ref{1.2}) is given by
$F_\phi(\eta):=\int_0^{2\pi}dx \mathrm{tr}\left(\phi(x)\Lambda(\eta(x))\right)$.
The formula of the Hamiltonian vector field of 
the admissible Hamiltonian $F_\phi$ on $\M_\CG$, 
 which is equivalent to (\ref{1.2}),
is described in \cite{BFP,Montreal}.
Because there is an underlying symplectic structure, there should exist 
a Poisson algebra containing the functions $F_\phi$ together with 
the functions of $M$ (and the Fourier coefficients of the differential polynomials 
in the affine Kac-Moody current $J= \eta' \eta^{-1}$).  
It is a non-trivial open problem to precisely characterize  such a Poisson algebra,
but fortunately this is not needed for the purposes of the present paper.
The key point to notice is that (\ref{1.15}) is well-defined as long as the matrix elements of $\Omega$
represent admissible Hamiltonians on $\M_\CG$, and this is the case 
since $\Omega$ depends only on the monodromy matrix $M$.
Note also that the problem of finding the momentum map in the finite dimensional context 
of the corresponding PL groupoid (with the PB in (\ref{2.22})) is literally the same as (\ref{1.15}) with 
the dependence on $x$ `forgotten'.

\section{The case of factorisable PL symmetry groups}
\setcounter{equation}{0}

We recall a convenient model of the Heisenberg double of a 
factorisable PL group in Subsection 2.1.
This will help us to recognize in Subsection 2.2 that the momentum map on
$\M_\CG$ is given by the monodromy matrix of the chiral WZNW field 
for the solution (\ref{1.13}), (\ref{1.14})  of the
PL-CDYBE (\ref{1.12}). We do this by converting $M$ into 
a `familiar variable' $\Omega$ on the dual PL group.
In Subsection 2.3 we use this result to show that the
two PL groupoids associated with
the solutions of (\ref{1.12}) in \cite{BFP} and in 
\cite{FM} are related by the same change of variables.

\subsection{Recall of the  Heisenberg double of factorisable PL groups}

Let $\G$ be a self-dual Lie algebra with a 
scalar product $\langle\ ,\ \rangle$.
We use the identification 
$\G\otimes \G \simeq \mathrm{End}(\G)$ defined by $\langle\ ,\ \rangle$
(whereby $X\otimes Y: Z \mapsto X \langle Y,Z\rangle$ for any $X,Y,Z\in \G$).
Suppose that $R^\nu\in \G\wedge\G$ ($\nu\neq 0$) is a solution of (\ref{1.7}) that has the form
\be
R^\nu = 2\nu R
\quad \hbox{with}\quad R\in \G\wedge \G.
\label{2.1}\ee
Recall that the dual group $G^*$ to $(G,\{\ ,\ \}_G^{R^\nu})$
is the subgroup of $G\times G$ corresponding to the Lie subalgebra
$\G^*$ of $\G\oplus \G$ given by
\be
\G^* = \{ (R^+(X), R^-(X))\,\vert\, X\in \G\}.
\label{2.2}\ee
Here $R^\pm = R\pm \frac{1}{2}I$, where $I$ is the identity operator on $\G$.
At the same time $\G$ is identified with the diagonal subalgebra
$\G\equiv \G^\delta =\{ (X,X)\vert X\in \G\} \subset \G\oplus \G$.
The subalgebras $\G^\delta$ and $\G^*$ are in duality with respect
to the following scalar product on $\G\oplus \G$:
\be
\langle\langle (X_1, Y_1), (X_2, Y_2) \rangle\rangle_{\nu} =
\frac{1}{2\nu} \left( \langle X_1,Y_1 \rangle - \langle X_2, Y_2\rangle \right).
\label{2.3}\ee
The normalization is chosen so that $\G\oplus \G$ with this scalar product
is the Drinfeld double of $\G$ with its Lie bialgebra structure
defined by $R^\nu$.

By definition, the Heisenberg double of $(G,\{\ ,\ \}_G^{R^\nu})$
is the space $G\times G^*$ equipped with the PB $\{\ ,\ \}^{R^\nu}$  given
in terms of the variables $g\in G$ and
$(\Omega^+, \Omega^-)\in G^*$ as follows:
\be
\{ g_1, g_2\}^{R^\nu}= 2\nu [g_1 g_2,R],
\label{2.4}\ee
\be
\{ \Omega^+_1, \Omega^-_2\}^{R^\nu} =2\nu [ R^+, \Omega^+_1 \Omega^-_2],
\quad
 \{ \Omega^\epsilon_1, \Omega^\epsilon_2\}^{R^\nu} = 2\nu
[R, \Omega^\epsilon_1 \Omega^\epsilon_2], \quad\epsilon=\pm,
\label{2.5}\ee
\be
\{ g_1, \Omega^\pm_2\}^{R^\nu} =-2\nu g_1 R^\mp \Omega^\pm_2.
\label{2.6}\ee
This PB is due to Semenov-Tian-Shansky \cite{STS}.
The Heisenberg double is
the natural PL analogue of the cotangent bundle $T^* G\simeq G \times \G^*$ 
(in the left trivialisation).
In fact \cite{STS}, (\ref{2.4}) and (\ref{2.5}) give the PBs on the PL groups
$G$ and $G^*$, respectively, while (\ref{2.6}) means (see below)
that $(\Omega^+,\Omega^-)$ serves as the momentum map for the natural
right PL action of $G$ on the phase space.

It is often convenient to identify $G^*$ with the open submanifold of $G$
provided by
\be
G_* := \{ \Omega = \Omega^+ (\Omega^-)^{-1} \,\vert\,
(\Omega^+, \Omega^-)\in G^*\}.
\label{2.7}\ee
The name factorisable PL group refers to the fact that the
elements of $G^*\simeq G_* \subset G$ are factorisable as above.
In terms of the variable $\Omega= \Omega^+ (\Omega^-)^{-1}$, the PBs in (\ref{2.5})
can be recast as the `Semenov-Tian-Shansky PB'
\be
\{ \Omega_1, \Omega_2\}^{R^v}=
 2\nu \left( R \Omega_1 \Omega_2 + \Omega_1 \Omega_2 R - \Omega_1 R^- \Omega_2 -
\Omega_2 R^+ \Omega_1\right).
\label{2.8}\ee
The PBs between $g$ and $(\Omega^+, \Omega^-)$ in (\ref{2.6}) are equivalent to
\be
\{ g_1, \Omega_2\}^{R^\nu} = 2\nu g_1( \Omega_2 R^+ - R^- \Omega_2).
\label{2.9}\ee
Formulae (\ref{2.8}), (\ref{2.9}) show that the PB on 
$G\times G_*$ smoothly extends to $G\times G$,
so that $G\times G$ can also be regarded as a PL analogue of the cotangent
bundle $G\times \G^*$.
We denote this latter version of the Heisenberg double by $(T^*G)_{R^\nu}$.

There is a natural right PL action of the group $(G, \{\ ,\ \}_G^{R^\nu})$ on
its Heisenberg double $((T^*G)_{R^\nu}, \{\ ,\ \}^{R^\nu})$. The action of any $q\in G$ is
given by the map $\R_q: (T^*G)_{R^\nu} \rightarrow (T^*G)_{R^\nu}$:
\be
\R_q: (g, \Omega) \mapsto (gq, q^{-1} \Omega q).
\label{2.10}\ee
Restricting $\Omega$ to $G_*$, let us now explain that
the non-Abelian momentum map for this action operates by mapping
$(g,\Omega)\in (T^*G)_{R^\nu}$ to its $\Omega$-component.
This obviously is a Poisson map.
It also generates the action,
since (as a consequence of (\ref{2.5}), (\ref{2.6}))
for any $T\in \G$ one can indeed express the infinitesimal
version of (\ref{2.10}) through the PB according to
\be
(gT)_{ij}= \langle\langle (T,T), \{ g_{ij},
 (\Omega^+,\Omega^-)\}^{R^\nu} (\Omega^+,\Omega^-)^{-1}\rangle\rangle_\nu
\label{2.11}\ee
and
\be
[\Omega, T]_{ij}=\langle\langle (T,T), \{ \Omega_{ij},
 (\Omega^+,\Omega^-)\}^{R^\nu} (\Omega^+,\Omega^-)^{-1}\rangle\rangle_\nu.
\label{2.12}\ee
The matrix elements used here refer to an arbitrary finite dimensional representation of $G$.

If $\Omega$ is near enough to $e\in G$, then we can uniquely
parametrize it as $\Omega = e^\omega$, where $\omega$ belongs to some
neighbourhood of zero in $\G$.
In the next subsection we will use the local expression of the PB on $(T^*G)_{R^\nu}$
in terms of this logarithmic variable.

\subsection{The momentum map on $\M_\CG$ in the factorisable case}

Our first main result is
the explicit description of the non-Abelian momentum map for the
PL symmetries on the chiral WZNW  phase space in the factorisable case,
which is given by the following theorem.

\medskip
\noindent
{\bf Theorem 1.}
{\em Consider the chiral WZNW phase space $(\M_\CG, \{\ ,\ \}_{WZ}^r)$ (\ref{1.1})
with the exchange $r$-matrix defined by eqs.~(\ref{1.10}), (\ref{1.13}), (\ref{1.14}).
Suppose that $R^\nu$ ($\nu\neq 0$) has the form (\ref{2.1})  and identify the
dual PL group $G^*$ with the domain $G_* \subset G$ (\ref{2.7})
equipped with the PB (\ref{2.8}).
Parametrize the monodromy matrix $M$ and
the variable $\Omega\in G_*$ as  $M=e^m$ and $\Omega =e^\omega$.
Then the non-Abelian momentum map $\M_\CG \rightarrow G_*$ associated with the
PL action in (\ref{1.9}) depends only on $M=\eta^{-1}(x)\eta(x+2\pi)$ and is
given explicitly  by
\be
M \mapsto \Omega(M)= M^{2\nu}
\quad \Longleftrightarrow \quad
m\mapsto \omega(m) = 2\nu m.
\label{2.13}\ee
Here it is assumed that $M$ lies near to $e\in G$
so that $\Omega(M)\in G_* \simeq G^*$. }

\medskip

In order to prove the theorem, we shall compare
the PBs on $\M_\CG$ and on $(T^*G)_{R^\nu}$ by using the
logarithmic variables $m$ and $\omega$.
The expressions to be compared are recorded in the next two lemmas.

\medskip
\noindent
{\bf Lemma 1.}
{\em Consider $(\M_\CG, \{\ ,\ \}_{WZ}^r)$
with the exchange $r$-matrix given by eqs.~(\ref{1.10}),
(\ref{1.13}), (\ref{1.14}).
In terms of $m:= \log M$, the PBs (\ref{1.3})
and (\ref{1.4}) can be equivalently written as
\be
\eta^{-1}(x) \{\eta(x),  \langle m, T\rangle \}_{WZ}^r=
\left(- R^\nu \circ \ad_{m}  + \chi(\nu\, \ad_{m}) \right)(T),
\label{2.14}\ee
\be
\{m, \langle m, T\rangle\}_{WZ}^r
= \left(- \ad_{m}  \circ R^\nu
+\chi(\nu\, \ad_{m}) \right)([m,T])
\label{2.15}\ee
for any $T\in \G$, where we use the analytic function $\chi$ defined in (\ref{1.14}).}

\medskip
\noindent
{\bf Lemma 2.}
{\em If $\Omega$ is near enough to $e\in G$, then the PBs (\ref{2.9})
and (\ref{2.8})
on the Heisenberg double $(T^* G)_{R^\nu}$ can be written in terms of
$\omega := \log \Omega$ and $\forall T\in \G$ as
\be
g^{-1} \{   g, \langle \omega, T\rangle \}^{R^\nu} = 2\nu
\left( -R \circ \ad_{\omega}  + \chi( \half \ad_{\omega})\right)(T),
\label{2.17}\ee
\be
\{ \omega, \langle \omega, T\rangle \}^{R^\nu}
= 2\nu \left( -\ad_{\omega}  \circ R
+\chi(\half \ad_{\omega})\right)([\omega,T]),
\label{2.18}\ee
with $\chi$ defined in (\ref{1.14}).}

\medskip
\noindent{\em Proof of Lemma 1 and Lemma 2.}
Consider a function $\psi\in C^\infty(G)$ and associate with it 
the function $\hat \psi$ on $\M_{\CG}$
by $\hat \psi: \eta \mapsto \psi(M)$, where 
$\eta(x+2\pi)=\eta(x)M$. From (\ref{1.3}) 
\be
\eta^{-1}(x)\{ \eta(x), \hat \psi\}^r_{WZ}= 
T^a (r^+_{ab}\cR^b \hat \psi - r^-_{ab} \cL^b \hat \psi),
\label{L1}\ee
where $\cR^b \hat \psi$ denotes the function on $\M_\CG$ as associated with
$\cR^b \psi\in C^\infty(G)$. The derivatives are defined in (\ref{1.6}) and $r^{\pm}_{ab}$ are
the matrix elements of the WZNW exchange $r$-matrix.
Similarly, for $\varphi,\psi\in C^\infty(G)$, (\ref{1.4}) gives
\be
\{ \hat \varphi, \hat \psi\}_{WZ}^r = r_{ab} (\cR^a \hat \varphi) (\cR^b \hat \psi) +
r_{ab} (\cL^a\hat \varphi)(\cL^b \hat \psi)
-r^-_{ab} (\cR^a\hat \varphi)(\cL^b \hat \psi)-r^+_{ab} (\cL^a\hat \varphi)(\cR^b \hat \psi).
\label{L2}\ee
Now define $\cL_X:= \langle X, T_a\rangle \cL^a$ and 
$\cR_X:=\langle X, T_a\rangle \cR^a$ for any constant $X\in \G$,
and consider the $\G$-valued function $M\mapsto \log M$ on 
an appropriate neighbourhood of $e\in G$.
Recall the following well-known (e.g. \cite{SW}) formulae:
\be
\cL_X \log M = \lambda(-\half \ad_{\log M})(X),
\qquad
\cR_X \log M= \lambda(\half \ad_{\log M})(X),
\label{L3}\ee
where $\lambda$ is the analytic function given by
\be
\lambda(z)={ze^z}(\sinh z)^{-1}.
\label{L4}\ee
Now it is a matter of direct calculation to 
obtain (\ref{2.14}) and (\ref{2.15}) from (\ref{L1}) and (\ref{L2})
by using (\ref{L3}) and the formula, (\ref{1.10}) with (\ref{1.13})-(\ref{1.14}), of $r(M)$.
The statement of Lemma 2 is verified in the same way. {\em Q.E.D.}

\bigskip
\noindent{\em Proof of Theorem 1.}
Consider the $\G$-valued function $\omega$ on $\M_\CG$ defined
in (\ref{2.13}). For this function,  the formulae (\ref{2.14})
and (\ref{2.15}) can be rewritten as
\be
\eta^{-1}(x) \{   \eta(x), \langle \omega, T\rangle \}_{WZ}^r = 2\nu
\left( -R \circ \ad_{\omega}  + \chi(\half \ad_{\omega})\right)(T),
\label{2.19}\ee
\be
\{ \omega, \langle \omega, T\rangle \}_{WZ}^r
= 2\nu \left(- \ad_{\omega}  \circ R
+\chi(\half \ad_{\omega})\right)([\omega,T]),
\qquad\forall T\in \G.
\label{2.20}\ee
Comparison of (\ref{2.18}) and (\ref{2.20}) shows that
$m\mapsto \omega(m) := 2 \nu m$ yields
a Poisson map from $(\M_\CG, \{\ ,\ \}_{WZ}^r)$ to
$(G^*, \{\ ,\ \}_{G^*}^{R^v})$ through the corresponding variables
$M=e^m$ and $\Omega=e^\omega$.
Since $\Omega$ serves as the momentum map for the PL action (\ref{2.10})
on $(T^* G)_{R^\nu}$, further comparison of (\ref{2.17}) with (\ref{2.19})
and (\ref{2.10}) with
(\ref{1.9}) shows that the map (\ref{2.13}) does indeed provide the required
non-Abelian momentum map on the chiral WZNW phase space.
Note that $\Omega(M)$ is guaranteed to lie in $G_*$ (\ref{2.7}) if $M$
is restricted\footnote{The map $M\mapsto \Omega(M)=M^{2\nu}$ would yield
a Poisson map to $G$ equipped with the Semenov-Tian-Shansky PB
(\ref{2.8}) even if $\Omega(M)$ lay outside the domain $G_* \simeq G^*$.}
to a suitable neighbourhood of $e\in G$.
{\em Q.E.D.}

\subsection{Connection between two PL groupoids}

Let us recall from \cite{BFP} that the WZNW exchange $r$-matrices
that appear in (\ref{1.2}) can be related to PBs on certain finite
dimensional Poisson manifolds as well.
Namely, on the manifold
\be
P_{\CG}:= \check G \times G \times \check G = \{ (\tilde M, g, M)
\,\vert\, M, \tilde M \in \CG,\, g\in G \},
\label{2.21}\ee
the following formula defines  a PB, $\{\ ,\ \}^{r}$,
for any solution of the $G$-CDYBE (\ref{1.5}):
\bea
&&
\{ g_1, g_2\}^r = g_1 g_2  r(M) -  r(\tilde M) g_1 g_2
\nonumber\\
&&  \{ g_1, M_2\}^r = g_1\bigl( M_2 r^+(M) - r^-(M) M_2\bigr)
\nonumber\\
&& \{ g_1, \tilde M_2\}^r =
\bigl( {\tilde M}_2 r^+(\tilde M) - r^-({\tilde M}) {\tilde M}_2\bigr) g_1
\nonumber\\
&& \{ M_1, M_2\}^r =
M_1M_2 r(M) + r(M)M_1 M_2 - M_1 r^-(M) M_2 - M_2 r^+(M) M_1
\nonumber\\
&& \{ \tilde M_1, \tilde M_2\}^r =
- \tilde M_1\tilde M_2 r(\tilde M) - r(\tilde M)\tilde M_1 \tilde M_2 +
\tilde M_1 r^-(\tilde M) \tilde M_2 + \tilde M_2 r^+(\tilde M) \tilde M_1
\nonumber\\
&& \{ M_1, \tilde M_2\}^r =0.
\label{2.22}
\eea
It can be shown \cite{BFP} that
-- with the  trivial groupoid structure --
$(P_\CG, \{\ ,\ \}^r)$ is a PL  groupoid in the sense of \cite{We}.
Note also in passing that if $r$ is associated with classical
$G$-symmetry on $\M_\CG$ as described around equation (\ref{1.16}),
then this PL groupoid coincides,
by means of the exponential parametrization of $M$ and $\tilde M$,
with the `dynamical PL groupoid'
of Etingof and Varchenko \cite{EV} that encodes the canonical $r$-matrix
obtained by taking
$\nu=0$ in (\ref{1.13})-(\ref{1.14}).

Motivated by the above and the geometric interpretation
of the CDYBE in \cite{EV}, we have proposed \cite{FM} a geometric
setting for the PL-CDYBE (\ref{1.12}) that relies on the Heisenberg double.
To present this, let us now denote the
elements of $P_{\hat G}$ as
\be
P_{\hat G}:= \hat G \times G \times \hat G =
\{ (\tilde \Omega, g, \Omega)\,\vert\, \Omega, \tilde \Omega \in \hat G,\, g\in G \},
\label{2.23}\ee
where $\hat G \subset G$ is some open submanifold.
Take a factorisable constant $r$-matrix,  $R^\nu=2\nu R$,  and
a  (smooth or holomorphic) map $\K: \hat G \mapsto \G\wedge \G$.
Then consider the following ansatz for a PB, $\{\ ,\ \}^{can}$, on $P_{\hat G}$:
\bea
&&\{ g_1, g_2\}^{can}= 2\nu \bigl( g_1 g_2 (R+\K(\Omega)) -
(R+\K(\tilde \Omega)) g_1 g_2\bigr)
\nonumber\\
&& \{ g_1, \Omega_2 \}^{can} =2\nu g_1 (\Omega_2 R^+ - R^- \Omega_2 )
\nonumber\\
&& \{ g_1, \tilde \Omega_2 \}^{can} =
2\nu  (\tilde \Omega_2 R^+ - R^- \tilde \Omega_2 )g_1
\nonumber\\
&& \{ \Omega_1, \Omega_2\}^{can}=2\nu\bigl(
 R \Omega_1 \Omega_2 + \Omega_1 \Omega_2 R
 -\Omega_1 R^- \Omega_2 -
\Omega_2 R^+ \Omega_1\bigr)
\nonumber\\
&& \{\tilde \Omega_1, \tilde \Omega_2\}^{can}=-2\nu\bigl(
R\tilde \Omega_1 \tilde \Omega_2 +
\tilde \Omega_1 \tilde \Omega_2 R  -
\tilde \Omega_1 R^- \tilde \Omega_2 - \tilde \Omega_2 R^+ \tilde \Omega_1\bigr)
\nonumber\\
&&\{ \Omega_1, \tilde \Omega_2\}^{can}=0.
\label{2.24}\eea
We assume that $\K$ is a {\em $G$-equivariant} map\footnote{The domains $\hat G$ and $\check G$ are 
chosen to be invariant under the adjoint action of $G$ on $G$.}, since this is
required locally around $e\in G$ by the Jacobi identity
$\{\{ g_1, g_2\}^{can}, \Omega_3\}^{can}
+ \hbox{cycl. perm}=0$ and its counterpart with $\tilde \Omega$.
Upon comparison with (\ref{2.8}) and (\ref{2.9}),
it is clear that (\ref{2.24}) defines a PB if the Jacobi identity
$\{\{g_1, g_2\}^{can}, g_3 \}^{can}+ \hbox{cycl. perm}=0$ holds.
This condition is found to be equivalent to the following
version of the PL-CDYBE:
\be
[R_{12}+\K_{12},R_{23}+ \K_{23}] +  T^a_1 (\half \D^+_{T_a} -\D_{R(T_a)}^-) \K_{23}
+\hbox{cycl. perm.} = {\I} \qquad\hbox{on}\qquad \hat G,
\label{2.25}\ee
where $\cal I$ is an arbitrary $G$-invariant constant element
of $\G\wedge \G\wedge \G$.
By using the equivariance of $\K$, the cross-terms containing both $\K$ and $R$ can
be cancelled.
If we set $\I=-\nu^2 \hat f$ and  $\K=-K^\nu$, then
(\ref{2.25}) becomes identical to equation (\ref{1.12}).

In the above setting the PL-CDYBE (\ref{2.25})
appears as the guarantee of the Jacobi
identity of the PB in (\ref{2.24}).
To further clarify the meaning of this interpretation,  first note that
$\Omega$ and $\tilde \Omega$ define in an obvious way the non-Abelian momentum maps
that generate natural PL actions of $(G, \{\ ,\ \}^{R^\nu}_G)$ on
$(P_{\hat G}, \{\ ,\ \}^{can})$ acting respectively
 by right- and left-multiplications on $g$:
\be
\R_q: (\tilde \Omega, g, \Omega) \mapsto (\tilde \Omega, g q, q^{-1} \Omega q)
\quad\hbox{and}\quad
{\bf L}_q: (\tilde \Omega , g, \Omega) \mapsto (q \tilde \Omega q^{-1}, q g, \Omega )
\quad \forall q\in G.
\label{2.26}\ee
The equivariance of $\K$ ensures that these are PL actions.
Second,  notice that the constraint $\tilde \Omega=g \Omega g^{-1}$
defines a Poisson 
submanifold of $(P_{\hat G}, \{\ ,\ \}^{can})$ for any dynamical $r$-matrix
$\K$, which is isomorphic with (an open submanifold of) the Heisenberg double
$(T^* G)_{R^\nu}$.
To put this differently, we may say that $(P_{\hat G}, \{\ ,\ \}^{can})$
is obtained from (an open submanifold of) the Heisenberg double by
`forgetting' the constraint  $\tilde \Omega=g \Omega g^{-1}$ between the
left- and right-momenta, and modifying the PB, by inserting $\K$
in the first line of (\ref{2.24}), in such a way to keep the
PL symmetries (\ref{2.26}).

\medskip
\noindent
{\em Remark 2.}  Another equivalent form of the PL-CDYBE (\ref{2.25}) is provided by
returning to the variable $(\Omega^+,\Omega^-)\in G^*$ and replacing $\K$ by
$\tilde \K: \check G^*\rightarrow \G\wedge \G$ defined as
\be
\tilde \K(\Omega^+, \Omega^-)=\K(\Omega)\quad \hbox{with}\quad \Omega=\Omega^+ (\Omega^-)^{-1},
\label{tildeK}\ee
where $\check G^*\subset G^*$ corresponds to $\hat G\subset G$.
Then, as is easy to check, the PL-CDYBE takes the more natural form
\be
[R_{12}+\tilde \K_{12}, R_{23}+\tilde \K_{23}] +  T^a_1 \cL_{T^*_a} \tilde \K_{23}
+\hbox{cycl. perm.} = {\I},
\label{PL-CDYBE}\ee
where $\cL_{T_a^*}$ denotes the derivative along the right-invariant vector field
on $G^*$ associated with the basis element $T_a^*= (R^+(T_a), R^-(T_a))\in \G^*$.
(The basis $\{T^a\}$ of $\G$ is dual to the basis $\{T_a^*\}$ of $\G^*$
with respect to the scalar product (\ref{2.3}) on the Drinfeld double with $\nu=1$.)
The use of the variable $(\Omega^+,\Omega^-)$ is conceptually more natural,
but in order to study the PL symmetries on the chiral WZNW phase space and {\em to find
solutions} of the PL-CDYBE, the use of the variable $\Omega$ seems more convenient.

\medskip

Let us now focus on the special cases of the PL groupoids $(P_\CG, \{\ ,\ \}^r)$
that belong to the WZNW exchange $r$-matrices associated with the
PL symmetry (\ref{1.9}).
In the finite dimensional setting (\ref{2.22}),  this symmetry
translates into right and left PL actions of
$(G, \{\ ,\ \}_G^{R^\nu})$ on $(P_\CG, \{\ ,\ \}^r)$
that operate quite similarly to (\ref{2.26}):
\be
\R_q: (\tilde M, g, M) \mapsto (\tilde M, gq, q^{-1} M q)
\quad\hbox{and}\quad
{\bf L}_q: (\tilde M , g, M) \mapsto (q \tilde M q^{-1}, qg, M )
\quad \forall q\in G.
\label{2.27}\ee
Theorem 1 implies the following proposition, which gives
the momentum maps for these PL actions and clarifies
the relationship between the PBs defined in (\ref{2.22}) and in (\ref{2.24}).

\medskip
\noindent
{\bf Proposition 1.}
{\em
Consider the Poisson manifolds $(P_{\CG},\{\ ,\ \}^r)$ and
$(P_{\hat G}, \{\ ,\ \}^{can})$ endowed with the PBs (\ref{2.22}) and (\ref{2.24})
defined respectively by
\be
r(M)= R^\nu + f_\nu(\ad_{\log M})\quad \hbox{and}\quad
\K(\Omega)= -f_{\frac{1}{4\nu}}(\ad_{\log \Omega}),
\label{2.28}\ee
where the function $f_\nu$ is given in (\ref{1.14}).
Here $\check G$ and $\hat G$ are open submanifolds
of $G$ for which the exponential parametrization is valid and
the respective functions $f_\nu(\ad_{\log M})$ and 
$f_{\frac{1}{4\nu}}(\ad_{\log \Omega})$
are well defined for $M\in \check G$ and $\Omega \in \hat G$.
We can choose these submanifolds 
in such a way that the map $M\mapsto M^{2\nu}$ 
yields a diffeomorphism from $\check G$ to $\hat G$.
Then the map
\be
(\tilde M,g, M)\mapsto (\tilde \Omega(\tilde M), g, \Omega(M))
:= (\tilde M^{2\nu}, g, M^{2\nu})
\label{2.29}\ee
is a Poisson diffeomorphism from $(P_{\CG},\{\ ,\ \}^r)$ to
$(P_{\hat G}, \{\ ,\ \}^{can})$.
Hence these PL groupoids are isomorphic. }

\medskip
\noindent
{\em Proof.}
As a consequence of Theorem 1,
the last five formulae of (\ref{2.22}) are translated into the
last five formulae of (\ref{2.24}) by the map (\ref{2.29}).
By expressing this map in the logarithmic variables,
\be
\omega(m)= 2\nu m,\,\,
\tilde \omega(\tilde m)= 2\nu \tilde m
\quad\hbox{with}\quad
m=\log M,\,\, \tilde m = \log \tilde M,\,\, \omega = \log \Omega,\,\,
\tilde \omega = \log\tilde \Omega,
\label{2.30}\ee
we can rewrite the first formula of (\ref{2.22}) as
\be
\{ g_1, g_2\}^r = 2\nu \bigl(g_1 g_2 (R - f_{\frac{1}{4\nu}}
(\ad_{\omega(m)}))-
(R - f_{\frac{1}{4\nu}}(\ad_{\tilde \omega(\tilde m)}))g_1 g_2\bigr),
\label{2.31}\ee
which agrees with the first formula of (\ref{2.24})
if $\K(\Omega)=-f_{\frac{1}{4\nu}}(\ad_\omega)$. {\em Q.E.D.}

\section{The case of the standard compact PL groups}
\setcounter{equation}{0}

Below we take $\G$ to be a compact simple Lie algebra equipped with
its standard, Drinfeld-Jimbo $r$-matrix, which is the most important
case from the point of view of physical applications of the WZNW model.
We next recall the relevant Heisenberg double, and then present the momentum map and some remarks
on the corresponding PL groupoids.
We concentrate on the logical outline of the statements, relegating
the underlying calculations to Appendix A.

\subsection{The Heisenberg double of a compact PL group}

Let $\A$ be a complex simple Lie algebra with a Chevalley basis 
given by $E_{\pm\alpha}$ ($\alpha\in \Phi^+$)
and $H_{\alpha_k}$ $(\alpha_k\in \Delta$), where $\Phi^+$ and $\Delta$ denote the set of positive
and simple roots, respectively.
With respect to  the Killing form $\langle\ ,\ \rangle$ of $\A$, normalized so that the long roots
have length $\sqrt{2}$,
one has $\langle E_\alpha, E_\beta\rangle = \frac{2}{\vert\alpha\vert^2} \delta_{\alpha, -\beta}$.
Let us now take $\G$ to be the compact real form of $\A$,
\be
\G={\mathrm{span}}_\bbR\{ \mathrm{i}(E_\alpha + E_{-\alpha}), 
(E_\alpha - E_{-\alpha}), \mathrm{i} H_{\alpha_k}\,\vert\,
\alpha\in \Phi^+,\,\alpha_k\in\Delta\}.
\label{3.1}\ee
Then the realification $\A_\bbR$ of $\A$ (i.e. $\A$ regarded as a 
Lie algebra over the reals) can be decomposed
as the vector space direct sum
\be
\A_\bbR= \G + \B
\label{3.2}\ee
with the `Borel subalgebra'
\be
\B={\mathrm{span}}_\bbR\{ E_\alpha, \mathrm{i} E_\alpha, H_{\alpha_k}
\,\vert\,
\alpha\in \Phi^+,\,\alpha_k\in\Delta\}.
\label{3.3}\ee
$\G$ and $\B$ are isotropic subalgebras with respect to the non-degenerate 
invariant bilinear form on $\A_\bbR$ defined by
\be
\langle\langle X,Y\rangle \rangle_\theta := \frac{1}{\theta} \mathrm{Im} \langle X, Y\rangle
\qquad \forall X,Y\in\A_\bbR\simeq \A,
\label{3.4}\ee
where $\theta\in \bbR$ is an arbitrary non-zero constant.
Thus $\A_\bbR$ carries the factorisable $r$-matrix
\be
\rho:=\frac{1}{2} ( \pi_\G - \pi_\B),
\label{3.5}\ee
 where $\pi_\G$ and $\pi_\B$ are the projections on
$\A_\bbR$ associated with the splitting (\ref{3.2}).
It is well known that the factorisable Lie bialgebra
$(\A_\bbR, \langle\langle\ ,\ \rangle\rangle_\theta, \rho)$
is the Drinfeld double of $\G$ equipped with its standard $r$-matrix,
given by $\Rit= \theta \Ri\in \G\wedge\G$ with
\be
\qquad
\Ri:= \sum_{\alpha \in \Phi^+} \frac{\vert \alpha\vert^2}{4}
(E_\alpha - E_{-\alpha})\wedge \mathrm{i}(E_\alpha + E_{-\alpha})=
\mathrm{i}\sum_{\alpha\in \Phi^+}\frac{\vert \alpha\vert^2}{2}  E_\alpha\wedge E_{-\alpha}.
\label{3.6}\ee
Our notation reflects the fact that $\Rit$ satisfies (\ref{1.7}) with $\nu=\mathrm{i}\theta$,
where $\hat f$ is defined by means of the restriction of the 
Killing form $\langle\ ,\ \rangle$ to $\G$.

Let $A_\bbR$ be a connected real Lie group with Lie algebra $\A_\bbR$ and denote
by
$G$ and $B$ the connected Lie subgroups associated with the subalgebras $\G$ and $\B$.
We equip the group $G$ with  the Sklyanin bracket $\{\ ,\ \}_G^{\Rit}$
written as
\be
\{ q_1, q_2\}_G^{\Rit}= \theta [q_1 q_2, \Ri].
\label{3.7}\ee
The dual PL group is $(B, \{\ ,\ \}_B^{\Rit})$, where the PB on $B$ is induced from
the Drinfeld double in the standard way.
The Heisenberg double of the compact PL group $(G,\{\ ,\ \}_G^{\Rit})$ is  the
Poisson (actually symplectic) space $(A_\bbR, \{\ ,\ \}^{\Rit})$, whose PB
can be described in the St Petersburg notation symbolically as follows \cite{STS}.
If $a$ denotes the $A_\bbR$-valued variable, then we have
\be
\{ a_1, a_2\}^{\Rit} =  -\hat \rho a_1 a_2 - a_1 a_2 \hat \rho,
\label{3.8}\ee
where $\hat \rho \in \A_\bbR\wedge \A_\bbR$ corresponds to $\rho\in \mathrm{End}(\A_\bbR)$ by means
of the scalar product $\langle\langle\ ,\ \rangle\rangle_\theta$.
This PB is further discussed in Appendix A.

We now use the Iwasawa and Cartan decompositions (e.g. \cite{OV}) 
of the group $A_\bbR$ to
produce a handy model of the Heisenberg double $(A_\bbR, \{\ ,\ \}^{\Rit})$.
By the Iwasawa decomposition, one can uniquely decompose any element $a\in A_\bbR$
according to
\be
a= g^{-1} \tilde b = b \tilde g\quad\hbox{with}\quad g, \tilde g \in G,\quad b,\tilde b \in B.
\label{3.9}\ee
As a manifold, we then identify $A_\bbR$ with $G\times B$ by the mapping $a\mapsto (g, b)$.
It can be shown that the map
$A_\bbR \rightarrow B$ that operates using (\ref{3.9}) as $a\mapsto b$  is the momentum map
for the right PL action of $(G,\{\ ,\ \}_G^{\Rit})$ on
$(A_\bbR,\{\ ,\ \}^{\Rit})$ defined by
\be
\R_q: a\mapsto q^{-1} a,
\qquad
\forall q\in G,\, a\in A_\bbR.
\label{3.10}\ee
Let us now trade the variable $b$ for the new variable
\be
\Omega:= b b^\dagger,
\label{3.11}\ee
where for $b=e^\beta$ we have $b^\dagger = e^{\beta^\dagger}$ with dagger standing
for minus one times the Cartan involution of $\A_\bbR$.
In other words, dagger is $-1$ on $\G$ and $+1$ on $\mathrm{i}\G$ in the Cartan decomposition
\be
\A_\bbR = \G + \mathrm{i}\G,
\label{3.12}\ee
which gives $E_\alpha^\dagger = E_{-\alpha}$, 
$H_{\alpha_k}^\dagger = H_{\alpha_k}$.
It follows from the corresponding Cartan decomposition\footnote{By the Cartan decomposition one can
uniquely write $a\in A_\bbR$ as $a=h\tilde p= p \tilde h$ with $h,\tilde h\in G$ and
$p,\tilde p\in e^{\mathrm{i}\G}$;
comparison with (\ref{3.9}) gives $aa^\dagger=b b^\dagger = p^2$, i.e., $p=e^{\mathrm{i}\omega}$.}
of the group $A_\bbR$
that we can uniquely parametrize $\Omega$ as
\be
\Omega= b b^\dagger = e^{2\mathrm{i} \omega} \quad\hbox{with}\quad \omega \in \G.
\label{3.13}\ee

Collecting the above mentioned facts, by using (\ref{3.9}) and 
(\ref{3.13}) we obtain a diffeomorphism
$A_\bbR \rightarrow G\times \G$ by the map
\be
a \mapsto (g, \omega)
\quad\hbox{with}\quad \omega = -\frac{\mathrm{i}}{2}\log \Omega,
\qquad
\Omega= aa^\dagger = b b^\dagger.
\label{3.14}\ee
Henceforth we use the pair $(g, \Omega)$ (or equivalently $(g,\omega)$)
as coordinates on $A_\bbR$.
In these convenient variables the PL action (\ref{3.10}) of $G$ on $A_\bbR$ takes the form
\be
\R_q: (g, \Omega) \mapsto (gq, q^{-1} \Omega q)
\qquad\forall q\in G,\, (g,\Omega) \in G\times e^{\mathrm{i}\G}\simeq A_\bbR.
\label{3.15}\ee
If we identify
$B$ and the domain $e^{\mathrm{i}\G}\subset A_\bbR$ as two models of $A_\bbR/G$,
given respectively by the Iwasawa and the Cartan
decompositions, then $\Omega$ yields directly the momentum map for this action.
The important point is that $b$ in (\ref{3.9}) can be uniquely recovered from
$\Omega= b b^\dagger$.

Now we are ready to present the key formula of this subsection.

\noindent
{\bf Lemma 3.}
{\em In terms of the coordinates $(g,\omega)\in G\times \G$ defined by (\ref{3.9}) and (\ref{3.13}),
the PB (\ref{3.8}) of the Heisenberg double $(A_\bbR,\{\ ,\ \}^{\Rit})$ takes the
following form:
\be
\{ g_1, g_2\}^{\Rit}= \theta [g_1 g_2, \Ri],
\label{3.16}\ee
\be
g^{-1} \{   g, \langle \omega, T\rangle \}^{\Rit} =\theta
\left( - \Ri\circ \ad_{\omega}  + \chi(\mi\, \ad_{\omega})\right)(T),
\label{3.17}\ee
\be
\{ \omega, \langle \omega, T\rangle \}^{\Rit}
= \theta \left(- \ad_{\omega}  \circ \Ri
+\chi(\mi\, \ad_{\omega})\right)([\omega, T]),
\label{3.18}\ee
where $T\in \G$ is an arbitrary constant 
and $\chi(\mi\, z) = z \cot z$. }
\medskip

Notice that the formulae in Lemma 3 are essentially the same as the ones in Lemma 2;
$\nu$ is purely imaginary in the compact case, $\omega$ is now 
defined by (\ref{3.13}), and $\chi(z)=z\coth z$ is replaced by
$\chi(\mi\, z)=z\cot z$.
The proof of Lemma 3 relies on a  routine but not quite trivial 
calculation and is sketched in Appendix A.

\subsection{Momentum map and PL groupoids in the compact case}

We describe here the analogues for the compact case
of the results of Subsections 2.2 and 2.3.
We can be very brief in what follows because, in terms of the convenient
variables introduced in Subsection 3.1, all
formulae involved are essentially the
same as in the factorisable case.

Since we know that the momentum map on $(A_\bbR,\{\ ,\ \}^\Rit)$ for the action (\ref{3.15})
is given by $(g, \Omega)\mapsto \Omega$, we obtain the momentum 
map for the PL symmetry (\ref{1.9}) on
$\M_{\check G}$ upon comparing the formulae 
in (\ref{2.14}), (\ref{2.15}) with those in (\ref{3.17}),
(\ref{3.18}).

\medskip
\noindent
{\bf Theorem 2.}
{\em Consider the chiral WZNW phase space $(\M_\CG, \{\ ,\ \}_{WZ}^r)$ (\ref{1.1})
with the exchange $r$-matrix defined by eqs.~(\ref{1.10}), (\ref{1.13}), (\ref{1.14}).
Suppose that $\nu=\mathrm{i}\theta$ with $\bbR\ni\theta\neq 0$ and
$R^\nu=\Rit=\theta \Ri$ with $\Ri$ in (\ref{3.6}).
Identify the dual, $G^*=B$, of the PL group $(G,\{\ ,\ \}_G^{\Rit})$ with the domain
$e^{\mathrm{i}\G}$ through the map 
$B\ni b\mapsto \Omega=bb^\dagger = e^{2\mathrm{i} \omega}\in e^{\mathrm{i}\G}$.
The non-Abelian momentum map $\M_\CG \rightarrow G^*$ associated with the
PL action in (\ref{1.9}) depends only on $M=\eta^{-1}(x)\eta(x+2\pi)$ and,
in terms of the parametrization  $M=e^m$, is given explicitly  by
\be
M \mapsto \Omega(M)= M^{2\nu} =e^{2\mathrm{i}\theta m}
\quad \Longleftrightarrow \quad
m\mapsto \omega(m) = \theta m.
\label{3.19}\ee
 }
\medskip
\noindent
{\em Proof.}
One readily verifies that by setting $\nu=\mathrm{i}\theta$, $R^\nu= \theta \Ri$ and $\omega=\theta m$,
the formulae in (\ref{2.14}), (\ref{2.15}) can be rewritten as
\be
\eta^{-1}(x) \{   \eta(x), \langle \omega, T\rangle \}_{WZ}^r =\theta
\left( - \Ri\circ \ad_{\omega}  + \chi(\mi\,\ad_{\omega})\right)(T),
\label{3.20}\ee
\be
\{ \omega, \langle \omega, T\rangle \}_{WZ}^r
= \theta \left(- \ad_{\omega}  \circ \Ri
+\chi(\mi\,\ad_{\omega})\right)([\omega, T]).
\label{3.21}\ee
These formulae have the same form as (\ref{3.17}) and (\ref{3.18}), respectively,
which obviously implies the statement of the theorem. {\em Q.E.D.}
\medskip

Motivated by the relation between the two PL groupoids described in Subsection 2.3,
we can convert the PL groupoid $(P_{\check G}, \{\ ,\ \}^r)$ associated by (\ref{2.22}) with the WZNW
exchange $r$-matrix into a canonical model of it in the compact case as well.
This  results by replacing the $\check G$-valued `monodromy variables' $M$ and $\tilde M$
by  corresponding $e^{\mathrm{i}\check \G}$-valued
`momentum variables'  $\Omega=(M)^{2\mathrm{i}\theta}$ and 
$\tilde \Omega=(\tilde M)^{2\mathrm{i}\theta}$.
Here $\check \G \subset \G$ is the open submanifold 
where $-\mathrm{i}\log\Omega$
takes its values, and we stress that 
$e^{\mathrm{i}\check \G}$ is a model of an open submanifold, $\check B$, of $B=G^*$.
Call the resulting PL groupoid $(P_{e^{\mathrm{i}\check \G}}, \{\ ,\ \}^{can})$ with
\be
P_{e^{\mathrm{i}\check \G}}:= e^{\mathrm{i}\check \G} \times G \times e^{\mathrm{i}\check \G} =
\{ (\tilde \Omega, g, \Omega)
\,\vert\, \Omega, \tilde \Omega \in e^{\mathrm{i}\check \G},\, g\in G \}.
\label{3.22}\ee
The PBs $\{g_1, \Omega_2\}^{can}$, $\{\Omega_1, \Omega_2\}^{can}$ and their `tilded variants'
are the same as for the Heisenberg double $(A_\bbR, \{\ ,\ \}^\Rit)$, and this also holds for
$\{\Omega_1, \tilde \Omega_2\}^{can}=0$ with the only difference being
that the relation $\tilde\Omega = g \Omega g^{-1}$ has now been `forgotten'.
By using these momentum variables the dynamical $r$-matrix only appears in the PB
\be
\{ g_1, g_2\}^{can}= \theta \bigl( g_1 g_2 (\Ri+\K(\Omega)) -
(\Ri+\K(\tilde \Omega)) g_1 g_2\bigr),
\label{3.23}\ee
where from (\ref{1.13}), (\ref{1.14}) now we have
\be
\K(\Omega)=\frac{1}{\theta} K^{\mathrm{i}\theta}(M)= \frac{1}{\theta}
f_{\mathrm{i}\theta}( (2\mi \theta)^{-1}\ad_{\log \Omega} )\in \mathrm{End}(\G)=\G\otimes \G,
\label{3.24}\ee
since $\log M= \frac{1}{2\mi \theta}\log\Omega$ by (\ref{3.19}).
In fact, taking (\ref{3.23}) as an ansatz with an unspecified function $\K$,
the Jacobi identities of the PB $\{\ ,\ \}^{can}$ turn out to be equivalent
to the following variant of the PL-CDYBE:
\be
[R^\mi_{12}+\K_{12},R^\mi_{23}+ \K_{23}] +  T^a_1 (\D_{\mi T_a}^+ - \D_{\Ri(T_a)}^-) \K_{23}
+\hbox{cycl. perm.} = {\I},
\label{3.25}\ee
where $\cal I$ is an arbitrary $G$-invariant constant element
of $\G\wedge \G\wedge \G$.
It is a good check on our arguments that $\K$ in (\ref{3.24}) indeed 
solves this equation,
with right-hand side $\I= (\frac{1}{(4\theta)^2}-\frac{3}{4})\hat f$.

Similarly to the factorisable case, the PL-CDYBE (\ref{3.25}) can again be rewritten in
a more natural form  by introducing $\tilde \K: \check G^*=\check B \rightarrow \G\wedge \G$ by the rule
\be
\tilde \K(b) = \K(\Omega)\quad\hbox{with}\quad \Omega=bb^\dagger\in e^{\mathrm{i} \check \G}.
\label{3.26}\ee
Let $\{ T_a^*\}\subset \B$ stand for the dual to a basis $\{ T^a\} \subset \G$, i.e.,
$\langle\langle T^a, T_b^*\rangle\rangle =\mathrm{Im} \langle T^a, T^*_b\rangle = \delta^{a}_{b}$.
As explained in Appendix A, in terms of $\tilde \K$ (\ref{3.25}) takes the form
\be
[R^\mi_{12}+\tilde \K_{12},R^\mi_{23}+ \tilde \K_{23}] +  T^a_1 \cL_{T^*_a} \tilde \K_{23}
+\hbox{cycl. perm.} = {\I},
\label{3.27}\ee
where $\cL_{T_a^*}$ denotes the derivative along the right-invariant vector field on $B$ associated
with $T_a^*\in \B$.
This equation is the same as (\ref{PL-CDYBE}).
It is clear that this version of the PL-CDYBE can be generalized to any PL group.

\section{Conclusion}
\setcounter{equation}{0}

In this paper we found the momentum map for the PL symmetry (\ref{1.9}) on the chiral
WZNW phase space $(\M_{\CG}, \{\ ,\ \}_{WZ}^r)$ equipped with the exchange $r$-matrix defined by
(\ref{1.10})-(\ref{1.14}).
The meaning of our result is that the momentum map is simply the monodromy matrix $M$
of the chiral WZNW field, in the sense that it is in one-to-one correspondence with a standard variable
$\Omega$ in
the dual of the PL symmetry group as  given explicitly by Theorem 1 and Theorem 2.
The correspondence between $M\in \check G\subset G$ and $\Omega\in \check G^*\subset G^*$
implies that the PL groupoids associated with the chiral 
WZNW phase 
space in \cite{BFP} are isomorphic
to what we call canonical PL groupoids,
which can be constructed starting from the Heisenberg double of 
$(G, \{\ ,\ \}_G^{R^\nu})$.
This  result is presented in Proposition 1 and around equation (\ref{3.23}).

Incidentally, a one-to-one correspondence between $M\in \check G$ and $\Omega \in \check G^*$ 
must clearly exist also in the cases that are not covered by
our analysis, for instance the cases with $\nu=0$ but $R^\nu\neq 0$ in (\ref{1.7}).
The comparison between (\ref{1.2})-(\ref{1.4}) and (\ref{2.22}) shows that  
finding the momentum map on the infinite dimensional manifold $(\M_\CG, \{\ ,\ \}^r_{WZ})$  
and on the finite dimensional manifold $(P_\CG, \{\ ,\ \}^r)$ are equivalent problems in these cases as well.

It turned out that in terms of the momentum variable $\Omega$ the exchange $r$-matrices that
were discovered in the WZNW model satisfy the following natural PL generalization of the CDYBE.
Let $\tilde \K: \check G^* \rightarrow \G\wedge \G$ be equivariant with respect to the pertinent
(dressing and adjoint) actions of the coboundary PL group $G$ with Lie algebra $\G$.
In the general case, the PL-CDYBE has the same form as (\ref{3.27}),
if $\cL_{T_a^*}$ denotes the derivative along the right-invariant vector field on $G^*$ associated
with $T_a^*\in \G^*$, where $T^a$ is a basis of $\G$, and $R^\mi$
is replaced by the element of $\G\wedge\G$ that defines the PL structure on $G$.
Clearly, the PL-CDYBE is the guarantee of the Jacobi identity of a PB on
a PL groupoid also in the general case, similarly to (\ref{2.24}).
Remarkably, this equation also arises in \cite{DM}, from some completely different considerations.
For self-dual Lie algebras the exchange $r$-matrices of \cite{BFP}, 
given by (\ref{1.10})-(\ref{1.14}),
always yield solutions of the PL-CDYBE after replacing the variable $M\in \check G$
by the variable $\Omega\in \check G^*$.
In fact, as will be detailed elsewhere, 
these solutions can be reduced to the duals of certain PL subgroups of $G$,
generalizing the Dirac reduction of solutions of the CDYBE studied 
in \cite{FGP}.
A different analysis of the reductions of (\ref{3.27}) is contained in \cite{EEM},
where the quantization of the PL dynamical $r$-matrices is also studied. 
It would be very interesting to apply 
the results of \cite{EEM} to develop the deformation quantization of the chiral WZNW phase space.

\bigskip
\bigskip
\noindent{\bf Acknowledgements.}
L.F. was supported in part by the Hungarian 
Scientific Research Fund (OTKA) under grant numbers 
T034170, T043159 and M036803.  We wish to thank T. Ratiu for hospitality at the EPFL.
We also wish to thank  J. Balog and B. Enriquez for useful comments on the manuscript.

\renewcommand{\theequation}{\arabic{section}.\arabic{equation}}
\renewcommand{\thesection}{\Alph{section}}
\setcounter{section}{0}

\section{Calculations on the `compact Heisenberg double'}
\setcounter{equation}{0}
\renewcommand{\theequation}{A.\arabic{equation}}

We here collect some useful, more or less well-known (see e.g.~\cite{LW}), formulae
concerning the Heisenberg double $( A_\bbR,\{\ ,\ \}^\Rit)$ used in Section 3. 

For any real function $\Phi\in C^\infty(A_\bbR)$ we define $D\Phi, D'\Phi \in C^\infty(A_\bbR,\A_\bbR)$ by
\be
\dt\Phi(e^{tX}ae^{tY})= 
\langle\langle (D\Phi)(a),X\rangle\rangle +\langle\langle (D'\Phi)(a), Y\rangle\rangle
\qquad  \forall X,Y\in \A_\bbR,\,a\in A_\bbR,
\label{A1}\ee
where $\langle\langle\ ,\ \rangle\rangle$ is given by (\ref{3.4}) with $\theta:=1$.
We often write $D_a\Phi$ for $(D\Phi)(a)$ and denote the adjoint action of $A_\bbR$ on $\A_\bbR$ simply by
conjugation.  
By the invariance of $\langle\langle\ ,\ \rangle\rangle$, we have  $D'_a \Phi = a^{-1} (D_a\Phi) a$.
Formula (\ref{3.8}) means that for $\Phi, \Psi\in C^\infty(A_\bbR)$
\be
\{\Phi,\Psi\}^\Rit(a)=-\theta \langle\langle D_a\Phi,\rho D_a\Psi\rangle\rangle -
\theta \langle\langle D'_a\Phi,\rho D'_a\Psi\rangle\rangle,
\label{A2}\ee
where $\rho$ is defined in (\ref{3.5}).
We next express this PB in terms of the coordinates $g\in G$ 
and $b\in B$ given by 
the factorisation 
(\ref{3.9}).
Identifying $\A_\bbR^*$ with $\A_\bbR$
by means of $\langle\langle\ ,\ \rangle\rangle$, and using (\ref{3.2}), 
$\G^*$ is naturally identified 
with $\B$ and $\B^*$ with $\G$. Hence, for functions
$\phi\in C^\infty(G)$ and $f\in C^\infty(B)$,
$D\phi, D'\phi\in C^\infty(G,\B)$ and $D f, D'f\in C^\infty(B,\G)$ are defined by
\bea
&&\dt\phi(e^{tX}ge^{tY})= \langle\langle D_g\phi,X\rangle\rangle
+\langle\langle D'_g\phi,Y\rangle\rangle
\qquad \forall X,Y\in\G,\, g\in G,\nonumber\\
&&\dt f(e^{tX}b e^{tY})= \langle\langle D_b f,X\rangle\rangle
+\langle\langle D'_b f,Y\rangle\rangle
\qquad \forall X,Y\in\B,\, b\in B.
\label{A3}\eea
With the projections that appear in (\ref{3.5}) we have
\be
D_g' \phi=\pi_\B\left( g^{-1} (D_g\phi) g\right),
\qquad 
D_b'f = \pi_\G\left( b^{-1} (D_b f) b\right),
\label{A.4}\ee
\be
D_g \phi=\pi_\B\left( g(D'_g\phi) g^{-1}\right),
\qquad 
D_b f = \pi_\G\left( b (D_b' f) b^{-1}\right).
\label{A4+}\ee
Referring to the Iwasawa decompositions (\ref{3.9}),
for any $\phi\in C^\infty(G)$ and $f\in C^\infty(B)$ we introduce $\hat \phi\in C^\infty(A_\bbR)$ and
$\hat f\in C^\infty(A_\bbR)$ by
\be
\hat \phi(a)= \phi(g)\quad\hbox{and}\quad
\hat f(a)=f(b) \quad\hbox{for}\quad a=g^{-1} \tilde b = b \tilde g\in A_\bbR.
\label{A5}\ee
It is then straightforward to check the following relations:
\be
D_a \hat \phi = - g^{-1} (D_g \phi) g,
\quad
D'_a \hat \phi = - \tilde b^{-1} (D_g \phi) \tilde b,
\quad
D_a \hat f = b (D'_b f) b^{-1},
\quad
D'_a \hat f = \tilde g^{-1} (D'_b f) \tilde g.
\label{A6}\ee
By using that $D_a'\hat\phi$ is $\B$-valued and $D_a'\hat f$ is $\G$-valued, it
is not difficult to obtain from (\ref{A2}) and (\ref{A6})  that
\be
\{ \hat \phi, \hat \psi\}^\Rit(a)= \theta \langle\langle D_g \phi, g( D'_g \psi) g^{-1}\rangle\rangle,
\qquad \forall \phi, \psi \in C^\infty(G),
\label{A7}\ee
\be
\{ \hat f, \hat h\}^\Rit(a)= \theta \langle\langle D_b f, b( D'_b h) b^{-1}\rangle\rangle,
\qquad \forall f, h \in C^\infty(B),
\label{A8}\ee
\be
\{ \hat \phi, \hat f\}^\Rit(a)=\theta \langle\langle D'_g \phi,  D_b f\rangle\rangle,
\qquad \forall \phi\in C^\infty(G),\, f \in C^\infty(B).
\label{A9}\ee
Equation (\ref{A7}) (resp. (\ref{A8})) means that the map
$A_\bbR \ni a\mapsto g\in G$ (resp. $A_\bbR \ni a\mapsto b\in B$) is a Poisson map from
$A_\bbR$ to $G$ (resp. to $B$) equipped with the PL structure appearing on the right hand side
of (\ref{A7}) (resp. (\ref{A8})).
Furthermore, one sees from (\ref{A9}) that $A_\bbR\ni a\mapsto b\in B$ is the momentum map
for the right PL action of $G$ on $A_\bbR$ defined in (\ref{3.10}).

Instead of the pair $(g,b)$,
we now rewrite the above PB in terms of the new coordinates $(g,\omega)\in G\times \G$
on $A_\bbR$, where  $\omega$ parametrizes $b\in B$ according to (\ref{3.13}).
For a function $\phi\in C^\infty(G)$, define $\nabla\phi, \nabla' \phi \in C^\infty(G,\G)$ by
\be
\dt\phi(e^{tX}ge^{tY})= \langle \nabla_g\phi ,X\rangle  + \langle \nabla'_g\phi, Y\rangle
\qquad  \forall X,Y\in \G.
\label{A10}\ee
On account of the invariance of $\langle\ ,\  \rangle$, we have $\nabla'_g\phi = g^{-1} (\nabla_g\phi) g$.

\medskip\noindent
{\bf Lemma 4.}
{\em Equation (\ref{A7}) can be equivalently expressed as
\be
\{\hat \phi,\hat \psi\}^\Rit(a)=\theta \langle \nabla'_g\phi,\Ri (\nabla'_g\psi)\rangle -
\theta \langle \nabla_g\phi,\Ri (\nabla_g\psi)\rangle,
\label{A11}\ee
where $\Ri\in \mathrm{End}(\G)$ is defined by (\ref{3.6}).
}

\smallskip
\noindent
{\em Proof.}
It follows from
\be
\langle \alpha^\dagger, \beta^\dagger\rangle =\overline{\langle \alpha,\beta \rangle},
\qquad\forall \alpha,\beta\in \A,
\label{A12}\ee
that
\be
\langle \langle X, Y\rangle\rangle = -\frac{1}{2} \langle X, \mathrm{i} (Y+Y^\dagger) \rangle,
\qquad
\forall X\in \G\,,Y\in \B.
\label{A13}\ee
This in turn implies
\be
\nabla_g \phi = -\frac{\mathrm{i}}{2}( D_g \phi + (D_g \phi)^\dagger ),
\qquad
\nabla'_g \phi = -\frac{\mathrm{i}}{2}( D'_g \phi + (D'_g \phi)^\dagger ).
\label{A14}\ee
Noting that any element of $\G$ can be uniquely represented in the form $\mathrm{i}(Y + Y^\dagger)$
with $Y\in \B$,  one readily verifies from formula (\ref{3.6}) that $\Ri\in \mathrm{End}(\G)$
operates according to
\be
\Ri: \mathrm{i} (Y+ Y^\dagger) \mapsto (Y^\dagger - Y),\qquad \forall Y\in \B.
\label{A15}\ee
As a consequence of (\ref{A15}),  relation (\ref{A14}) is equivalent to
\be
D_g \phi = \mathrm{i} \nabla_g \phi +\Ri(\nabla_g \phi),
\qquad
D'_g \phi = \mathrm{i} \nabla'_g \phi +\Ri(\nabla'_g \phi).
\label{A16}\ee
Now for any $\alpha, \beta\in \B$ and $g\in G$ it is easy to check (from (\ref{A12}) using also
$g^\dagger = g^{-1}$) that
\be
4\langle\langle \alpha, g\beta g^{-1}\rangle\rangle =
 \langle \mathrm{i}(\beta +\beta^\dagger), g^{-1} (\alpha^\dagger -\alpha) g\rangle
+\langle \mathrm{i}(\alpha +\alpha^\dagger), g (\beta^\dagger -\beta) g^{-1}\rangle .
\label{A17}\ee
Equation (\ref{A11}) results by applying (\ref{A17}) to $\alpha=D_g\phi$, $\beta=D'_g\psi$ on the
right hand side of (\ref{A7}) taking (\ref{A16}) and the antisymmetry of $\Ri$ into account.
{\em Q.E.D.}
\smallskip

Clearly, Lemma 4 is equivalent to equation (\ref{3.16}). In order to obtain equations (\ref{3.17}) and
(\ref{3.18}), we associate with any $f\in C^\infty(B)$ the function $\tilde f\in C^\infty(\G)$ by
\be
f(b)= \tilde f(\omega ) \quad\hbox{with}\quad b b^\dagger = e^{2\mi \omega},
\label{A18}\ee
and express (\ref{A8}) and (\ref{A9}) in terms of $\tilde f$ and $\tilde h$.
Any $\tilde f\in C^\infty(\G)$  has the  $\G$-valued gradient $d\tilde f$  defined by
\be
\dt \tilde f(\omega + tX )= \langle d_\omega \tilde f ,X\rangle
\qquad  \forall X\in \G.
\label{A19}\ee
We can then verify the following statement.

\medskip\noindent
{\bf Lemma 5.}
{\em For $f\in C^\infty(B)$ and $\tilde f\in C^\infty(\G)$ related by (\ref{A18})
\be
D_b f = \left(-\Ri\circ \ad_\omega  + \chi(\mi\,\ad_\omega)\right)(d_\omega \tilde f),
\label{A20}\ee
\be
b (D'_b f) b^{-1}= \mi [\omega, d_\omega \tilde f] +  \chi(\mi\, \ad_\omega)(d_\omega \tilde f),
\label{A21}\ee
with the function $\chi(\mi z) = z \cot z$ as introduced in (\ref{1.14}).
}

\smallskip\noindent
{\em Proof.}
We start by noting that
\be
-2\langle \langle D_b f, Y\rangle\rangle =\langle D_b f, \mi (Y+Y^\dagger)\rangle =
\langle d_\omega \tilde f, \mi  \lambda(-\ad_{\mi \omega})(Y) +
\mi \lambda(\ad_{\mi \omega})(Y^\dagger)\rangle,
\quad \forall Y\in \B,
\label{A22}\ee
where the function $\lambda$ is given in (\ref{L4}).
The second relation is a consequence of (\ref{A18}) and
$(e^{tY} b) (e^{tY} b)^\dagger = e^{tY} (b b^\dagger) e^{t Y^\dagger}$.
Proceeding as in the proof of Lemma 4, one then shows that
\be
\mi  \lambda(-\ad_{\mi \omega})(Y) +\mi \lambda(\ad_{\mi \omega})(Y^\dagger)= -
[\omega, R^\mi(X)] + \chi(\ad_{\mi\omega})(X),
\qquad X=\mi(Y+Y^\dagger)\in \G.
\label{A23}\ee
These relations imply (\ref{A20}).
The proof of (\ref{A21}) is rather similar, and we omit it.
{\em Q.E.D.}
\smallskip

Incidentally, (\ref{A20}) and (\ref{A21}) are consistent on account of the identity
\be
\Ri(X)= \pi_\G(-\mi X)\qquad \forall X\in \G.
\label{A24}\ee
Clearly, (\ref{A20}) and (\ref{A21}) imply the following result.

\medskip\noindent
{\bf Lemma 6.}
{\em Using the above notations, the PBs in (\ref{A8}) and (\ref{A9}) can be rewritten as
\be
\{ \hat \phi, \hat f\}^\Rit(a)= \theta \langle \nabla'_g \phi,
\left(-\Ri\circ\ad_\omega  + \chi(\mi\,\ad_\omega)\right)(d_\omega \tilde f)\rangle,
\label{A25}\ee
\be
\{ \hat f, \hat h\}^\Rit(a)= \theta \langle d_\omega \tilde f,
\left(-\ad_\omega \circ \Ri  + \chi(\mi\,\ad_\omega)\right)([\omega, d_\omega \tilde h])\rangle.
\label{A26}\ee
}

The statements given in Lemmas 4 and 6 are equivalent to the formulae in Lemma 3,
which we used in Section 3.

Finally, let us explain the equivalence between equations (\ref{3.25}) and (\ref{3.27}).
For this,  consider a pair of functions $\F\in C^\infty(e^{\mi \G})$ and $\tilde\F\in C^\infty(B)$
related by
\be
\tilde\F(b) = \F(\Omega) \quad\hbox{for}\quad \Omega=b b^\dagger,\,\, b\in B.
\label{A27}\ee
For any $X\in \G$, the derivative that appears in (\ref{3.25}) is given by
\be
(\D_{\mi X}^+ \F- \D_{\Ri(X)}^- \F )(\Omega)= \dt \F(\Omega_t^X),
\qquad
\Omega^X_t = e^{t(\mi X +\Ri(X))} \Omega e^{t(\mi X -\Ri(X))},
\label{A28}\ee
which is well defined since $\Omega_t^X\in e^{\mi \G}$ for any real $t$.
One can check the identity
\be
\Omega_t^X = b_t^{-2Y} (b_t^{-2Y })^\dagger
\quad\hbox{with}\quad  b_t^{-2Y} = e^{-2 Yt} b,
\quad X=\mi(Y+ Y^\dagger),\, Y\in \B.
\label{A29}\ee
This implies that
\be
(\D_{\mi X}^+ \F- \D_{\Ri(X)}^- \F )(\Omega) = (\cL_{-2 Y} \tilde \F)(b),
\label{A30}\ee
where $\cL_{-2 Y} \tilde \F \in C^\infty(B)$ is the natural left-derivative of $\tilde \F$.
Now $T^b\in \G$ and $T_a= i(Y_a + Y_a^\dagger)\in \G$ satisfy
$\langle T^b,T_a\rangle = \langle\langle T^b, -2 Y_a\rangle\rangle$ by (\ref{A13}), which means
that $T_a^*=-2Y_a$ in the notation appearing in (\ref{3.27}).
Equation (\ref{3.27}) arises from the above consideration by taking $\F$ to be
the $\G\wedge\G$ valued function $\K$.


\begin{thebibliography}{EnGH}

\bibitem{Wi}
E. Witten, 
{\sl Non-abelian bosonization in two dimensions},
CMP  {\bf 92} (1984),  455-472.

\bibitem{CFT} P. Di Francesco, P. Mathieu and D. S\'en\'echal,
Conformal Field Theory (Springer, 1997).
 
\bibitem{Fad} L. Faddeev, {\sl On the exchange matrix of the WZNW model},
CMP {\bf 132} (1990), 131-138.

\bibitem{BDF}
J. Balog, L. Dabrowski and L. Feh\'{e}r, 
{\sl Classical $r$-matrix and exchange algebra in WZNW and Toda theories},
Phys. Lett. B {\bf 244} (1990), 227-234. 

\bibitem{AF}
A. Alekseev  and L. Faddeev, 
{\sl $(T^*G)_t$: a toy model for conformal field theory},
CMP {\bf 141} (1991), 413-422.

\bibitem{FG}
F. Falceto and K. Gawedzki, 
{\sl Lattice Wess-Zumino-Witten model and quantum groups},
J. Geom. Phys. {\bf 11} (1993), 251-279 
(arXiv hep-th/9209076).

\bibitem{Feld}
G. Felder, {\sl Conformal field theory and integrable systems 
associated with elliptic curves},
in: {Proc. of the ICM 94}, Birkhauser, 1994, pp.~1247-1255
(arXiv: hep-th/9407154).

\bibitem{EV}  P. Etingof and A. Varchenko, 
{\sl Geometry and classification of solutions of the classical 
dynamical Yang-Baxter equation},
CMP {\bf 192} (1998), 77-129 (arXiv: q-alg/9703040).

\bibitem{ES}
P. Etingof and  O. Schiffmann, 
{\sl Lectures on the dynamical Yang-Baxter equations},
arXiv: math.QA/9908064.

\bibitem{E}
P. Etingof, {\sl On the dynamical Yang-Baxter equation},
arXiv: math.QA/0207008.

\bibitem{BFP}
J. Balog, L.  Feh\'er  and  L. Palla, 
{\sl The chiral WZNW phase space and its Poisson-Lie groupoid},
Phys. Lett. B  {\bf 463} (1999), 83-92 (arXiv: hep-th/9907050);\\
J. Balog, L.  Feh\'er and  L. Palla, 
{\sl Chiral extensions of the WZNW phase space, 
Poisson-Lie symmetries and groupoids},
Nucl. Phys. B {\bf 568} (2000), 503-542 (arXiv: hep-th/9910046).

\bibitem{DM}
J. Donin and A. Mudrov,
{\sl Dynamical Yang-Baxter equation and quantum vector bundles},
arXiv: math.QA/0306028;\\
J. Donin and A. Mudrov,
{\sl Quantum groupoids and dynamical categories},
arXiv: math.QA/0311316.

\bibitem{Montreal}
J. Balog, L. Feher and L. Palla,
{\sl On the chiral WZNW phase space, exchange $r$-matrices and Poisson-Lie groupoids}, 
in: CRM Proceedings and Lectures Notes, 
Volume 26, eds.~J. Harnad et al, AMS, 2000, pp.~1-19 
(arXiv: hep-th/9912173).

    
\bibitem{BFPinPLA}
J. Balog, L. Feh\'er and L. Palla, 
{\sl The chiral WZNW phase space as a quasi-Poisson space},
Phys. Lett. A {\bf 277} (2000), 107-114
(arXiv:hep-th/0007045).

\bibitem{AKM}
A. Alekseev,Y. Kosmann-Schwarzbach and  E. Meinrenken, 
{\sl  Quasi-Poisson manifolds},
Canad. J. Math. {\bf 54} (2002), 3-29
(arXiv:math.DG/0006168). 

\bibitem{FM}
L. Feh\'er and I. Marshall,
{\sl On a Poisson-Lie analogue of the classical 
dynamical Yang-Baxter equation for self-dual Lie algebras},  
Lett. Math. Phys. {\bf 62} (2002), 51-62 
(arXiv:math.QA/0208159).

\bibitem{Lu}
J.-H. Lu, {\sl Momentum mappings and reduction of Poisson actions},
pp. 209-226 in: Symplectic Geometry, Groupoids, and Integrable Systems
(Berkeley, 1989), MSRI Publ., vol. 20 (Springer, 1991). 

\bibitem{BB} 
O. Babelon and D. Bernard, 
{\sl Dressing symmetries}, 
CMP {\bf 149} (1992), 279-306 (arXiv: hep-th/9111036).

\bibitem{SoibAMS}
L.I. Korogodski and Y.S.  Soibelman, 
Algebras of Functions on Quantum Groups: Part I (AMS, 1998). 

\bibitem{AM} 
A. Alekseev and  E. Meinrenken, 
{\sl The non-commutative Weil algebra},
Invent. Math. {\bf 139} (2000), 135-172 (arXiv: math.DG/9903052).

\bibitem{PF}
B.G. Pusztai  and L. Feh\'er, 
{\sl A note on a canonical dynamical $r$-matrix},
J. Phys. A  {\bf 34} (2001), 10949-10962 (arXiv: math.QA/0109082).

\bibitem{STS}
M.A. Semenov-Tian-Shansky, 
{\sl Dressing transformations and Poisson-Lie group actions}, 
Publ. RIMS {\bf 21} (1985), 1237-1260;\\
M.A. Semenov-Tian-Shansky, 
{\sl Poisson-Lie groups, quantum duality principle and the quantum double}, 
Theor. Math. Phys. {\bf 93} (1992), 
1292-1307 (arXiv: hep-th/9304042).

\bibitem{BD}
A.A. Belavin and V.G. Drinfeld,
{\sl Solutions of the classical Yang-Baxter equation for simple Lie
algebras}, Funct. Anal. Appl. {\bf 16} (1982), 159-180. 

\bibitem{Rawns}
M. Cahen, S. Gutt and J. Rawnsley, 
{\sl Some remarks on the classification of Poisson Lie groups},
Contemp. Math. {\bf 179} (1994) 1-16.

\bibitem{Soib}
Y.S. Soibelman, {\sl Algebra of functions on a compact quantum group
and its representations}, Leningrad Math. J. {\bf 2} (1991), 193-225. 

\bibitem{SW}
D.H. Sattinger and  O.L. Weaver,  
Lie Groups and Algebras with Applications to Physics, Geometry, 
and Mechanics (Springer, 1986).

\bibitem{We}
A. Weinstein, 
{\sl Coisotropic calculus and Poisson groupoids},  
J. Math. Soc. Japan {\bf 40} (1988), 705-727.

\bibitem{OV}
A.L. Onischik and E.B. Vinberg, 
Lie Groups and Lie Algebras III, Encyclopaedia of Mathematical 
Sciences, Vol. 41 (Springer, 1994).

\bibitem{FGP}
L. Feh\'er,  A. G\'abor and  B.G. Pusztai, 
{\sl On dynamical $r$-matrices obtained from Dirac reduction 
and their generalizations to affine Lie algebras},
J. Phys. A {\bf 34} (2001), 7235-7248 (arXiv: math-ph/0105047). 

\bibitem{EEM}
B. Enriquez, P. Etingof and I. Marshall,
{\sl Quantization of some Poisson-Lie dynamical $r$-matrices and Poisson homogeneous spaces},
arXiv: math.QA/0403283.


\bibitem{LW}
J.-H. Lu and A. Weinstein,
{\em Poisson Lie groups, dressing transformations and Bruhat decompositions},
J. Diff. Geom. {\bf 31} (1991), 501-526.


\end{thebibliography}
\end{document}